\documentclass[10pt,leqno]{amsart}
\usepackage{indentfirst, amssymb,amsmath,amsthm}
\usepackage{csquotes}
\usepackage{graphicx}
\usepackage{epstopdf}
\graphicspath{{G:\Back up\My Folders\My Work\Papers\Capsule Notes\Pythagorian theorem}}
\newcommand{\be}{\begin{equation}}
\newcommand{\ee}{\end{equation}}
\newcommand{\beas}{\begin{eqnarray*}}
\newcommand{\eeas}{\end{eqnarray*}}
\newcommand{\bea}{\begin{eqnarray}}
\newcommand{\eea}{\end{eqnarray}}

\numberwithin{equation}{section}
\AtBeginDocument{{\noindent\small  Resonance 21, 653–655 (2016)\\
DOI: 10.1007/s12045-016-0369-6}}
\begin{document}
\title[Pythagorean Theorem from Heron's Formula-Another Proof]{ Pythagorean Theorem from Heron's Formula-Another Proof }
\date{}
\author[B. Chakraborty ]{Bikash Chakraborty }
\date{}
\address{ Department of Mathematics, University of Kalyani, West Bengal 741235, India.}
\email{bikashchakraborty.math@yahoo.com, bikashchakrabortyy@gmail.com}
\maketitle
\let\thefootnote\relax
%%%%%%%%%%%%%%%%%%%%%%%%%%%%%%%%%%%%%%%%%%%%%%%%%%%%%%%%%%%%%%%%%%%%%%%%%%%%%%%%%%%%%%%%%%%%%%%%%%%%%%%%%%%%%%%
\begin{abstract} In this article using elementary school level Geometry  we observe an alternative proof of Pythagorean Theorem from Heron's Formula.
\end{abstract}
%%%%%%%%%%%%%%%%%%%%%%%%%%%%%%%%%%%%%%%%%%%%%%%%%%%%%%%%%%%%%%%%%%%%%%%%%%%%%%%%%%%%%%%%%%%%%%%%%%%%%%%%%%%%%%%%%%%%%%%%
\section{ Pythagorean Theorem from Heron's Formula -Another Proof }
For a right angled triangle ABC (where $\angle BAC$ is the right angle), we have to show that $$\overline{AB}^{2}+\overline{AC}^{2}=\overline{BC}^{2}.$$
%%%%%%%%%%%%%%%%%%%%%%%%%%%%%%%%%%%%%%%%%%%%%%%%%%%%%%%%%%%%%%%%%%%%%%%%%%%%%%%%%%%%%%%%%%%%%%%%%%%%%%%%%%%%%%%%%%%%%%%%%%%%%%%%%%%%%%%%%%%%%%%%%%%%%%%%%%
%\begin{figure}[h]
%\includegraphics{CIRCLE}
%\centering
%\end{figure}
\begin{figure}[h]
\centering
\includegraphics[scale=.7]{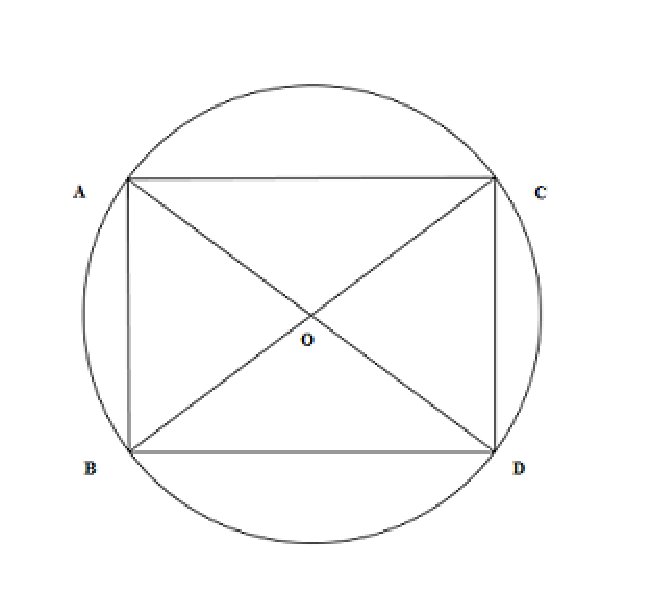}
\caption{}
\end{figure}
%%%%%%%%%%%%%%%%%%%%%%%%%%%%%%%%%%%%%%%%%%%%%%%%%%%%%%%%%%%%%%%%%%%%%%%%%%%%%%%%%%%%%%%%%%%%%%%%%%%%%%%%%%%%%%%%%%%%%%%%%%%%%%%%%%%%%%%%%%%%%%%%%%%%%%%%%%%%
\\
For the proof, let O be the mid point $\overline{BC}$. Since $\angle BAC$ is the right angle, we can draw a circle with center O and radius $\overline{BO}$ such that the points $A,B,C$ lies on the circle. Now we join two points $A$ and $O$ and draw a line segment from the point B parallel to AC which meets the circle at D, (say). Again we join the points C,D and O,D.\par
Assume $\overline{AB}=a,~\overline{AC}=b,~\overline{BO}=r$. \\
Clearly from the figure, \bea\label{1} &&The~ area~ of~ the~ rectangle~ABDC\\
\nonumber &=& 2(Area~ of~ the ~triangle~ ABO ~+~Area~ of~ the~ triangle~ AOC).\eea
Now the area of the triangle ABO =$|\sqrt{s(s-a)(s-r)(s-r)}|=|\frac{a}{2}\sqrt{r^{2}-\frac{a^{2}}{4}}|,$\\
where $s=\frac{a+r+r}{2}$, and $|x|$ denotes the absolute value of $x$.\\
Similarly the area of the triangle AOC= $|\frac{b}{2}\sqrt{r^{2}-\frac{b^{2}}{4}}|.$\\
Thus from \ref{1},
$$ab=|a\sqrt{r^{2}-\frac{a^{2}}{4}}|+|b\sqrt{r^{2}-\frac{b^{2}}{4}}|.$$
By putting $t=4r^{2}-a^{2}-b^{2}$ in above equation we have
\bea\label{2} && 2ab=|a\sqrt{t+b^{2}}|+|b\sqrt{t+a^{2}}|,\eea
i.e.,
\beas  2a^{2}b^{2}-(a^{2}+b^{2})t=2ab|\sqrt{t^{2}+(a^{2}+b^{2})t+a^{2}b^{2}}|
\eeas
Again squaring both sides, we get
\beas
(a^{2}-b^{2})^{2}t^{2}-8a^{2}b^{2}(a^{2}+b^{2})t=0.
\eeas
Now we consider two cases :\\
\textbf{Case-1} $a=b$\\
Then clearly $t=0$, i.e., $4r^{2}=a^{2}+b^{2}$.\\
\textbf{Case-2} $a\not=b$\\
In this case we have $t\geq 0$.\\
Now when $t> 0$, then
$$|a\sqrt{t+b^{2}}|+|b\sqrt{t+a^{2}}|>2ab,$$

which contradicts the equation (\ref{2}).\\
So $t=0$, i.e., $4r^{2}=a^{2}+b^{2}$.\\
Hence $\overline{AB}^{2}+\overline{AC}^{2}=a^{2}+b^{2}=(2r)^{2}=\overline{BC}^{2}$.
Hence the proof.

\end{document}